\def\half{{\textstyle{1\over 2}}}

\def\squareforqed{\hbox{\rlap{$\sqcap$}$\sqcup$}}
\def\qed{\ifmmode\squareforqed\else{\unskip\nobreak\hfil
\penalty50\hskip1em\null\nobreak\hfil\squareforqed
\parfillskip=0pt\finalhyphendemerits=0\endgraf}\fi\medskip}
\def\supp{{\rm supp\,}}
\def\bbbf{{\rm I\!F}}
\def\bbbc{{\mathchoice
{\setbox0=\hbox{$\displaystyle\rm C$}\hbox{\hbox
to0pt{\kern0.4\wd0\vrule height0.9\ht0\hss}\box0}}
{\setbox0=\hbox{$\textstyle\rm C$}\hbox{\hbox to0pt{\kern0.4\wd0\vrule
height0.9\ht0\hss}\box0}} {\setbox0=\hbox{$\scriptstyle\rm
C$}\hbox{\hbox to0pt{\kern0.4\wd0\vrule height0.9\ht0\hss}\box0}}
{\setbox0=\hbox{$\scriptscriptstyle\rm C$}\hbox{\hbox
to0pt{\kern0.4\wd0\vrule height0.9\ht0\hss}\box0}}}}
\def\bbbz{{\mathchoice {\hbox{$\sf\textstyle Z\kern-0.4em Z$}}
{\hbox{$\sf\textstyle Z\kern-0.4em Z$}} {\hbox{$\sf\scriptstyle
Z\kern-0.3em Z$}} {\hbox{$\sf\scriptscriptstyle Z\kern-0.2em Z$}}}}

\documentstyle[12pt]{article}
 
\newtheorem{theorem}{Theorem} 
\newtheorem{lemma}{Lemma} 
 
\begin{document}
\title{ The ratio and generating function of cogrowth
coefficients of finitely generated
groups} \author{Ryszard Szwarc\thanks{\noindent This work has been
partially supported by KBN (Poland) under grant 2 P03A 030 09.}}
\date{}
 \maketitle
\begin{abstract}
Let G be a group generated by $r$ elements $g_1,g_2,\ldots, g_r.$ Among
the reduced words in $g_1,g_2,\ldots, g_r$ of length $n$ some, say
$\gamma_n,$ represent the identity element of the group $G.$ It has been
shown in a combinatorial way
that the $2n$th root of $\gamma_{2n}$ has a limit, called the cogrowth
exponent  with respect to generators $g_1,g_2,\ldots, g_r.$ We show
by analytic methods that the numbers $\gamma_n$ vary regularly; i.e.
the ratio $\gamma_{2n+2}/\gamma_{2n}$ is also
convergent.  Moreover we derive new precise information on the domain of
holomorphy of $\gamma(z),$ the generating
function associated with the coefficients $\gamma_n.$
\end{abstract}
\footnotetext{\noindent 1991 {\it Mathematics Subject Classification}.
Primary  20F05, 20E05}
\footnotetext{\noindent
{\it Key words and phrases}: cogrowth of subgroups, free group, amenable
groups}

Every group $G$ generated by $r$ elements can be realized as a
quotient of the free group $\bbbf_r$ on $r $ generators by a normal
subgroup $N$ of $\bbbf_r,$ in such a way that the generators of the
free group $\bbbf_r$ are sent to the generators of the group $G.$
With the set of generators of $\bbbf_r$ we associate the length function
of words in these generators.  The cogrowth coefficients $\gamma_n=
\#\{x\in N \mid |x|=n\} $
were first introduced by Grigorchuk in \cite{g}.  The numbers $\gamma_n$
measure how big the group $G$ is when compared with $\bbbf_r.$
It has been shown that the quantities $\sqrt[2n]{\gamma_{2n} }$ have a
limit denoted by $\gamma,$ and called the growth exponent of $N$ in
$\bbbf_r.$ Since the subgroup $N$ can have at most $2r(2r-1)^{n-1}$
elements
of length $n,$ the cogrowth exponent $\gamma$ can be at most $2r-1.$ The
famous Grigorchuk result,
proved independently by J. M. Cohen in \cite{c}, states that the group
$G$ is amenable if and only if $\gamma=2r-1$ (see also \cite{sz1},
\cite{wo}).

The main result of this note is that the coefficients $\gamma_{2n}$
satisfy not only the Cauchy $n$th root test but also the d'Alambert
ratio test.
\begin{theorem}
The ratio of two consecutive even cogrowth coefficients
 $\gamma_{2n+2}/\gamma_{2n}$ has a limit.  Thus the
ratio tends to $\gamma^2,$ the square of the cogrowth exponent.
 \end{theorem}

{\it Proof.}
Let us denote by $g_1,g_2,\ldots,g_r$ the generators of $G.$ Let $\mu$
be the  measure equidistributed over the generators and
their inverses according to the formula
$$\mu = {1\over 2\sqrt{q}}\sum_{i=1}^{r}(g_i+g_i^{-1}),$$
where $q=2r-1.$ By an easy transformation of \cite[Formula (*)]{sz1} we
obtain
\begin{equation}\label{imp}
{z\over 1-z^2}\sum_{n=0}^{\infty}\gamma_nz^n
= {1\over 2\sqrt q}\sum_{n=0}^{\infty}\mu^{*n}(e) \left ({2\sqrt q
z\over qz^2+1}\right )^{n+1} ,
\end{equation}
for small values of $|z|.$
Let $\varrho$ denote the spectral radius of the random walk defined
by $\mu;$ i.e.
$$\varrho = \lim_{n\to\infty}\sqrt[2n]{\mu^{*2n}(e)}.$$
By $d\sigma(x)$ we will denote the spectral measure of this random walk.
Hence
\begin{equation}\label{sigma}
\mu^{*n}(e) =\int_{-\varrho}^\varrho x^n d\sigma(x).
\end{equation}
Note that the point $\varrho$ belongs to the support of $\sigma.$
Combining (\ref{imp}) and (\ref{sigma}) gives
\begin{eqnarray}
{z\over 1-z^2}\sum_{n=0}^{\infty}\gamma_nz^n&=&
{1\over 2\sqrt q}\int_{-\varrho}^\varrho \sum_{n=0}^{\infty}x^n \left
({2\sqrt q z\over
qz^2+1}\right )^{n+1}d\sigma(x)        \nonumber\\
&=&{1\over 2\sqrt q}
\int_{-\varrho}^{\varrho}   {z \over 1-2\sqrt{q}
xz+qz^2}d\sigma(x).\label{wazne}
\end{eqnarray}
By the well known formula for the generating function of the second kind
Chebyshev polynomials $U_n(x)$ (see \cite[(4.7.23), page 82]{s}) where
\begin{equation}\label{cz}
U_n(\half (t+t^{-1}))= {t^{n+1}-t^{-n-1}\over t-t^{-1}},
\end{equation}
we have
$$ {1\over 1-2\sqrt q xz+qz^2}= \sum_{n=0}^{\infty}U_n(x)q^{n/2}z^{n}.$$
Thus
$${z\over 1-z^2}\sum_{n=0}^{\infty}\gamma_nz^n= z
\sum_{n=0}^{\infty}q^{n/2}z^{n}\int_{-\varrho}^\varrho U_n(x)
d\sigma(x) .  $$
Therefore for $n\ge 2$ we have
\begin{equation}\label{najw}
\gamma_n=
q^{n/2} \int_{-\varrho}^\varrho \{U_n(x)- q^{-1}U_{n-2}(x)\}\,d\sigma(x).
\end{equation}
Since $U_{2n}(-x)=U_{2n}(x)$ we get
\begin{equation}\label{naj}
\gamma_{2n}=
q^{n}   \int_0^\varrho \{U_{2n}(x)-
q^{-1}U_{2n-2}(x)\}\,d\widetilde{\sigma}(x),
\end{equation}
where $\widetilde{\sigma}(A)=\sigma(A)+\sigma(-A)$ for $A\subset
(0,\varrho]$ and $\widetilde{\sigma}(\{0\})=\sigma(\{0\}).$
Let $$
I_n=\int_0^\varrho \{U_{2n}(x)-
q^{-1}U_{2n-2}(x)\}\,d\widetilde{\sigma}(x) .$$
By \cite[Corollary 2]{k} we have
$\varrho > 1.$
Hence we can split the integral $I_n$
into two integrals: the first $I_{n,1}$ over the interval
$[0,\varrho_0]$
and the second $I_{n,2}$ over
$[\varrho_0, \varrho],$ where $\varrho_0=(1+\varrho)/2.$
By (\ref{cz}) we have $|U_m(x)|\le (m+1)$ for $x\in[0,1]$ and
$$|U_m(x)|\le (m+1)[x+\sqrt{x^2-1}]^m \quad \mbox{for }x\ge 1.$$ Thus we
get
\begin{eqnarray}
I_{n,1} &\le &
2(2n+1)\left (\varrho_0+ \sqrt{\varrho_0^2-1}\right
)^{2n}\int_0^{\varrho_0}d\widetilde{\sigma}(x)\nonumber
\\
& \le   & 2(2n+1)\left
(\varrho_0+ \sqrt{\varrho_0^2-1}\right
)^{2n}
.  \label{est1}
\end{eqnarray}

Let's turn to estimating the integral $I_{n,2}$ over $[\varrho_0,\varrho].$
By (\ref{cz}) one can easily check that
$$\left | U_n(x)- {(x+\sqrt{x^2-1})^{n+1}\over 2\sqrt{x^2-1}}\right |
= o(1) \quad\mbox{when} \ n\to \infty,$$
uniformly on the interval $[\varrho_0,\varrho].$
Hence
$$ \left | U_{2n}(x)-q^{-1}U_{2n-2}(x)- (x+\sqrt{x^2-1})^{2n-1}
{(x+\sqrt{x^2-1})^{2}-q^{-1}\over
2\sqrt{x^2-1}}\right | =o(1),$$
when $n$ tends to infinity, uniformly in the interval
$[\varrho_0,\varrho].$
This implies
\begin{equation}\label{est2}
I_{n,2} \approx \widetilde{I}_{n,2}=
\int_{\varrho_0}^\varrho (x+\sqrt{x^2-1})^{2n}{(x+\sqrt{x^2-1})^{2}-q^{-1}\over
2\sqrt{x^2-1}(x+\sqrt{x^2-1})}\,d\widetilde\sigma(x).
\end{equation}
Since the endpoint $\varrho$ belongs  to the support of
$\widetilde{\sigma},$ we get
\begin{equation}\label{root}
{\widetilde{I}_{n,2}}^{1/{2n}} \longrightarrow \varrho +
\sqrt{\varrho^2-1}.
\end{equation}
By combining this with (\ref{est1}) and (\ref{est2}) we obtain
\begin{equation}\label{appr}
I_n= I_{n,1}+I_{n,2}= \widetilde{I}_{n,2}(1+o(1)),\qquad n\to
\infty.
\end{equation}
In view of   (\ref{root}) the integral $\widetilde{I}_{n,2}$ tends to
infinity. Thus by (\ref{naj}) and (\ref{appr})
 we have
$${\gamma_{2n+2}\over \gamma_{2n}} \approx
q{\widetilde{I}_{n+1,2}\over \widetilde{I}_{n,2}} .$$

\begin{lemma}[\cite{w}]\label{lem}
Let $f(x)$ be a positive and continuous function on $[a,b],$ and $\mu$ be
a finite measure on $[a,b].$ Then
$$\lim_{n\to \infty}{\int_a^bf(x)^{n+1}d\mu(x)\over
\int_a^bf(x)^{n}d\mu(x)}  =\max \{f(x)\mid\, x\in \supp \mu \}.$$
\end{lemma}

Applying Lemma 1 and using the fact that $\varrho$ belongs to the
support of $\widetilde{\sigma}$ gives
\begin{equation}\label{last}
{\gamma_{2n+2}\over \gamma_{2n}}\to
q\left \{\varrho+\sqrt{\varrho^2-1}\right \}^2.\qquad \qquad
\end{equation}
\qed
\begin{theorem}
The generating function $\gamma(z)=\sum_{n=0}^{\infty}\gamma_nz^n$ can
be decomposed into a sum of two functions $\gamma^{(0)}(z)$ and
$\gamma^{(1)}(z)$ such that $\gamma^{(0)}(z)$ is analytic in the  open disc of
radius $q^{-1/2}$  (where $q=2r-1),$ while $\gamma^{(1)}(z)$ is
analytic in
the whole complex plane after removing the two real intervals
$[-\gamma q^{-1},\, -\gamma^{-1}]$
and $[\gamma^{-1}, \gamma q^{-1}].$ Moreover,
$\gamma^{(1)}$ satisfies the functional equation
$${z\gamma^{(1)}(z)\over 1-z^2}= {(q/z)\gamma^{(1)}(q/z)\over (q/z)}.$$
\end{theorem}

{\it Proof.}
By (\ref{wazne}) we have
$$\gamma(z)=(1-z^2)\int_{-\varrho}^{\varrho}    {1 \over 1-2\sqrt{q}
xz+qz^2}d\sigma(x).$$
Let
\begin{eqnarray*}
    \gamma^{(0)}(z)&=& (1-z^2)\int_{-1}^{1} {1 \over 1-2\sqrt{q}
xz+qz^2}d\sigma(x),\\
 \gamma^{(1)}(z)&=& (1-z^2)\int_{1<|x|\le \varrho}  {1 \over 1-2\sqrt{q}
xz+qz^2}\,d\sigma(x).
\end{eqnarray*}
For $-1\le x\le 1$ the expression
$1-2\sqrt{q} xz+qz^2$ vanishes only on the circle of radius
$q^{-1/2}. $ Thus $\gamma^{(0)}(z)$ has the desired property.
For $1< |x|\le \varrho$ the expression
$1-2\sqrt{q} xz+qz^2$ vanishes only on the intervals
$$\left [-{\varrho+\sqrt{\varrho^2-1}\over \sqrt{q}},\,
-{\varrho-\sqrt{\varrho^2-1}\over \sqrt{q}}\right ],\quad  \left  [
{\varrho-\sqrt{\varrho^2-1}\over \sqrt{q}},\,
{\varrho+\sqrt{\varrho^2-1}\over \sqrt{q}}\right ].$$
By (\ref{last}) we have that
$\gamma=q^{1/2}(\varrho+\sqrt{\varrho^2-1}).$
This shows that $\gamma^{(1)}$ is analytic where it has been required.

The functional equation follows immediately from the formula
$${z\gamma^{(1)}(z)\over 1-z^2}=
\int_{1<|x|\le \varrho} \ {1 \over z^{-1}-2\sqrt{q}
x+qz}\,d\sigma(x).$$
\qed
\noindent{\bf Remark.} Combining (\ref{naj}) and (\ref{appr}) yields
$$\gamma_{2n}= q^n\left \{\int_{\varrho_0}^\varrho
(x+\sqrt{x^2-1})^{2n}{(x+\sqrt{x^2-1})^{2}-q^{-1}\over
2\sqrt{x^2-1}(x+\sqrt{x^2-1})}\,d\widetilde\sigma(x)+ o(1) \right \}.$$
We have
\begin{eqnarray*}
 h(\varrho_0):= {(\varrho_0+\sqrt{\varrho_0^2-1})^{2}-q^{-1}\over
2\sqrt{\varrho_0^2-1}(\varrho_0+\sqrt{\varrho_0^2-1})}&\ge &
{(x+\sqrt{x^2-1})^{2}-q^{-1}\over
2\sqrt{x^2-1}(x+\sqrt{x^2-1})} ,\\
  {(\varrho+\sqrt{\varrho^2-1})\over \varrho}\,x
&\ge &  x+\sqrt{x^2-1}.
\end{eqnarray*}
Therefore, in view of (\ref{sigma}), we get
\begin{eqnarray*}
\gamma_{2n} &\le & q^n \left \{
h(\varrho_0){(\varrho+\sqrt{\varrho^2-1})^{2n}\over
\varrho^{2n}}
\int_0^\varrho x^{2n}d\widetilde{\sigma}(x)+o(1)\right \}\\
& = & q^n
h(\varrho_0)\left
\{
{(\varrho+\sqrt{\varrho^2-1})^{2n}}
{\mu^{*2n}(e)\over \varrho^{2n}}+o(1)\right \}.
\end{eqnarray*}
Finally we obtain
$${\gamma_{2n}\over \gamma^{2n}}{\varrho^{2n}\over \mu^{*2n}(e)}=
{\gamma_{2n}\over \mu^{*2n}(e)}\left \{{\varrho\over\sqrt q (\varrho
+\sqrt{\varrho^2-1}) }\right \}^{2n} \le h(\varrho_0)+o(1).$$
We conjecture that the opposite estimate also holds; i.e.  the quantity
on the left hand side is bounded away from zero, by a positive constant
depending only on $\varrho.$ This conjecture can be checked easily if
the measure $\sigma$ is smooth in the neighbourhood of
$\varrho$ and the density has zero of finite order at $\varrho.$
\bigskip

{\bf Acknowledgement.}
For a long time I thought Theorem 1 follows from the following
statement.
\begin{verse}
{\em If $f(z)=\sum_{n=0}^{\infty}a_nz^n$ is analytic in the complex
plane  except\\ the half line $[1,+\infty),$ then the ratio $a_{n+1}/a_n$
converges to $1.$}
\end{verse}
I am grateful to Jacek Zienkiewicz from my Department for constructing a
fine
counterexample to this statement.

\bigskip

{\baselineskip 12pt
 \begin{flushleft}
Institute of Mathematics\\
Wroc\l aw University\\
pl. Grunwaldzki 2/4\\
50--384 Wroc\l aw, Poland\\
{\tt e-mail: szwarc@math.uni.wroc.pl}
\end{flushleft}}

\end{document}